\documentclass[reqno]{amsart}

\usepackage[utf8]{inputenc}

\usepackage[figuresright]{rotating}

\usepackage{amssymb}
\usepackage{amsthm}

\usepackage{setspace}
\newcommand{\R}{\mathbb{R}}
\newcommand{\C}{\mathbb{C}}
\newcommand{\f}{\rightarrow}
\newcommand{\deb}{\bar\partial}
\newcommand{\de}{\partial}
\newcommand{\K}{K\"{a}hler}
\newcommand{\ngh}{neighbourhood}

\newcommand{\lmb}{\lambda}
\newcommand{\ov}[1]{\overline{#1}}

\newcommand{\D}{\mathcal{D}}

\newcommand{\Aut}{\operatorname{Aut}}

\newcommand{\Ent}{\operatorname{Ent_d}}

\newcommand{\Isom}{\operatorname{Isom}}

\newcommand{\ol}{\operatorname{Hol}}

\newcommand{\ep}{{\varepsilon}}

\newcommand{\hilb}{\mathcal{H}}

\newcommand{\HBD}{M}

\newcommand{\gs}{\mathfrak{s}}
\newcommand{\g}{h}
\newcommand{\set}[2]{ \left\{#1\,:\,#2\right\} }

\newtheorem{prop}{Proposition}
\newtheorem{thm}{Theorem}
\newtheorem{lem}[prop]{Lemma}
\newtheorem{cor}[prop]{Corollary}
\newtheorem{defn}[prop]{Definition}

\newtheorem{ex}[prop]{Example}

\begin{document}

\author[R. Mossa]{Roberto Mossa}

\address{Dipartimento di Matematica e Informatica, Palazzo delle Scienze, via Ospedale 72 - 09124 Cagliari}

\email{roberto.mossa@gmail.com}

\title[A note on diastatic entropy and balanced metrics]
{A note on diastatic entropy and balanced metrics}

\begin{abstract}
We give un upper bound $\Ent\left( \Omega,\, g\right)<\lambda$  of the diastatic entropy $\Ent\left( \Omega,\, g\right)$ (defined by the author in \cite{mossa1}) of a complex bounded domain $\left( \Omega,\, g\right)$ in terms of the balanced condition (in Donaldson terminology) of the \K\  metric $\lambda\, g$.   When $\left( \Omega,\, g\right)$ is a homogeneous bounded domain
we   show that the converse holds true, namely if  $\Ent\left( \Omega,\, g\right)<1$ then $g$ is balanced. Moreover,  we explicitly compute   $\Ent(\Omega, g)$ in terms of Piatetski-Shapiro constants.
\end{abstract}

\maketitle


\section{Introduction and statements of the main results}
Let $\left( M, \, g \right)$ be a real analytic \K\ manifold with associated \K\ form $\omega$. Fix a local coordinate system $z=(z_1, \dots, z_n)$ on a \ngh\ $U$ of a point $p\in M$ and let $\phi: U \f \R$ be a real analytic \K\ potential (i.e $\omega_{|U}= \frac{i}{2}\, \de\deb\, \phi$). By shrinking $U$, we can assume that the potential $\phi\left( z\right) $ can be analytically continued to $U \times U$. Denote this extension by $\phi \left( z, \, \ov w\right)$. The Calabi's \emph{diastasis function} $\D:U\times U \f \R$ (E. Calabi \cite{calabi}) is defined by 
\[
\D\left(z,\, w  \right)= \D_w\left( z\right) = \phi \left( z, \, \ov z\right) + \phi \left( w, \, \ov w\right) - \phi \left( z, \, \ov w\right) - \phi \left( w, \, \ov z\right).
\]
Assume that $\omega$ has global diastasis function $\D_p:M \f \R$ centered at $p$. The \emph{diastatic entropy} at $p$ is defined as
\begin{equation}\label{def dent}
\Ent\left( M, \, g\right)\left( p\right)= \min \left\{ c >0 \ | \, \int_M e^{-c \, \D_p }\, \frac{\omega^n}{n!}<\infty \right\} .
\end{equation}
This definition does not depend on the point $p$ fixed, provided that for every $p\in \Omega$ the diastasis $\D_p$ is globally defined (see \cite[Proposition 2.2]{mossa1}).
The concept of diastatic entropy was defined by the author in \cite{mossa1} (following the ideas developed in \cite{exponential}),  where he obtained upper and lower bound for the first eigenvalue of a real analytic \K\ manifold.

In this paper we study the link between the diastatic entropy $\Ent \left( \Omega,\, g\right)  $ and the balanced condition of the metric $g$. S. Donaldson in \cite{donaldson}, in order to obtain a link between the constant scalar curvature condition on a \K\ metric $g$ and the Chow stability of a polarization $L$, gave the definition of a balanced metric on a compact manifold. Later, this definition has been generalized, by C. Arezzo and A. Loi in \cite{arlcomm} (see also \cite{arlquant}), to the noncompact case. 

In this paper we are interested in complex domains $\Omega\subset\C^n$ (connected open subset of $\C^n$) endowed with a real analytic K\"ahler metric $g$. Let $\omega$ be the \K\ form associated to $g$.
Assume that $\omega$ admits a global \K\ potential $\phi :\Omega\rightarrow\R$. We can define the weighted Hilbert space $\hilb_\phi$ of square integrable holomorphic functions on $\Omega$ with weight $e^{-\phi}$
\begin{equation}\label{hilbertspacePhi}
\hilb_\phi=\left\{ s\in\ol(\Omega)  :  \int_\Omega e^{-\phi}\left|s\right|^2\,\frac{\omega^n}{n!}<\infty\right\}.
\end{equation}
If $\hilb_\phi\neq \{0\}$ we can pick an orthonormal basis $\{s_j\}_{j=1, \dots, N \leq \infty},$
\[
\int_\Omega e^{-\phi}\,s_j \,\ov s_k\, \frac{\omega^n}{n!}=\delta_{jk}
\]
and define its reproducing kernel by
\begin{equation}\label{eq kernel}
K_{\phi}\left( z, \ov w\right) =\sum_{j=0}^{N} s_j\left( z\right)\, \ov {s_j\left( w\right) } .
\end{equation}
The so called \emph{$\ep$-function} is defined by
\begin{equation}\label{epsilon}
\varepsilon_g\left( z\right) =e^{-\phi\left(  z\right) }\,K_{\phi}\left( z, \,\bar z\right),
\end{equation}
as suggested by the notation $\varepsilon _g$ depends only on the metric $g$ and not on the choice of the K\"ahler potential $\phi$ (see, for example,  \cite[Lemma 1]{lm1} for a proof).
\begin{defn}\label{def bal}\rm
The metric $g$ is \emph{balanced} if and only if the function $\varepsilon_g$ is a positive  constant.
\end {defn}
\noindent In the literature the constancy of $\varepsilon_g$ was studied first by   J.  H. Rawnsley in \cite{rawnsley} under the name of $\eta$-{\em function}, later renamed as $\varepsilon$-{\em function} in \cite{cgr1}. It  also appears under the name of {\em distortion function} for the study of Abelian varieties by  J. R. Kempf \cite{kempf} and S. Ji \cite{ji}, and for complex projective varieties by S. Zhang \cite{zhang}. (See also  \cite{lm2} and references therein).

The following theorem represents the first result of this paper.
\begin{thm}\label{prop ent bal}
Let $\Omega\subset\C^n$ be a complex domain and $g$ a \K\ metric such that $\lmb\, g$ is balanced for some $\lmb >0$. Then the diastatic entropy $\Ent\left (\Omega,\,  g \right)$ is constant and $$\Ent\left (\Omega,\,  g \right)<\lmb.$$
\end{thm}

In particular (see Lemma \ref{ultimo lemma} below) the diastatic entropy of a homogeneous bounded  domain is constant.

It is natural to ask when the converse of the previous theorem holds true. In the next theorem we show that this is the case  when $\left( \Omega,\, g\right)$ is assumed to be homogeneous.
\begin{thm} \label{thm dentropy000} 
Let $\left( \Omega,\, g\right) $ be a homogeneous bounded  domain, then $g$ is balanced if and only if $$\Ent\left (\Omega,\,  g \right)<1.$$
\end{thm}
\noindent
We also compute the diastatic entropy of a homogeneous bounded domain in terms of the constants
$p_k,$ $b_k,$ $q_k$ and $\gamma_k$ defined by the Piatetski-Shapiro's root structure (see \eqref{def const} and \eqref{def const2} in Appendix below): 

\begin{cor}\label{cor1}
Let $\left( \Omega,\, g \right)$ be a homogeneous bounded domain. Then the diastatic entropy  of $\left( \Omega , \, g \right) $ is the positive constant given by
\begin{equation*} \label{eqn:def_of_lambda0}
 \Ent\left(\Omega,\, g\right)\left( z \right)  = \max_{1 \le k \le r} \frac{1 + p_k + b_k + q_k/2}{\gamma_k},\qquad \forall \, z \in \Omega.
\end{equation*}
\end{cor}

\begin{ex}\rm
If $g$ is the Bergman metric on $\Omega$,
 then $\gamma_k = 2 + p_k + q_k + b_k,$ $k=1, \dots, r$
 (see \cite[Theorem 5.1]{G64} or \cite[(2.19)]{N01}),
 so by Theorem  \ref{cor1},
 \begin{equation}\label{eq berg ent}
 \Ent\left(\Omega,\, g\right) \left( z \right) = \max_{1 \le k \le r} \frac{1 + p_k + b_k + q_k/2}{2 + p_k + q_k + b_k}, \qquad \forall \, z \in \Omega.
 \end{equation}
In particular,
 when $\left( \Omega,\, g\right)$ is a bounded symmetric domain (i.e. is a homogeneous bounded domain and a symmetric space as Riemannian manifold), there exist integers $a$ and $b$
 such that
 \[
 p_k = \left( k-1\right) a, \quad
 q_k = \left( r-k\right) a, \quad
 b_k = b, \quad
 \gamma_k = \left( r-1\right) a + b + 2.
 \]
where $r$ is the rank of $\left(  \Omega,\, g\right)$. Therefore, denoted by $\gamma$ the genus of $\Omega$, by \eqref{eq berg ent} we obtain
\begin{equation*}
\Ent\left(\Omega,\, g\right)\left( z \right) = \max_{1 \le k \le r} \frac{1 + \left( k-1\right) a + b + \left( r-k\right)a/2}{\gamma}
\end{equation*}
\begin{equation*}
 = \frac{1 + \left( r-1\right) a + b}{\gamma} = \frac{\gamma - 1}{\gamma}, \qquad \forall \, z \in \Omega.
\end{equation*}
See also \cite{mossa2} for a similar formula relating the volume entropy with the invariants $a$, $b$, and $r$.
\noindent The table below summarizes the numerical invariants and the dimension of irreducible bounded symmetric domains according to their type (for a more detailed description of these invariants, see e.g. \cite{arazy}, \cite{zhang}).

\begin{figure}[htb!]
\begin{center}
\includegraphics[width=126mm]{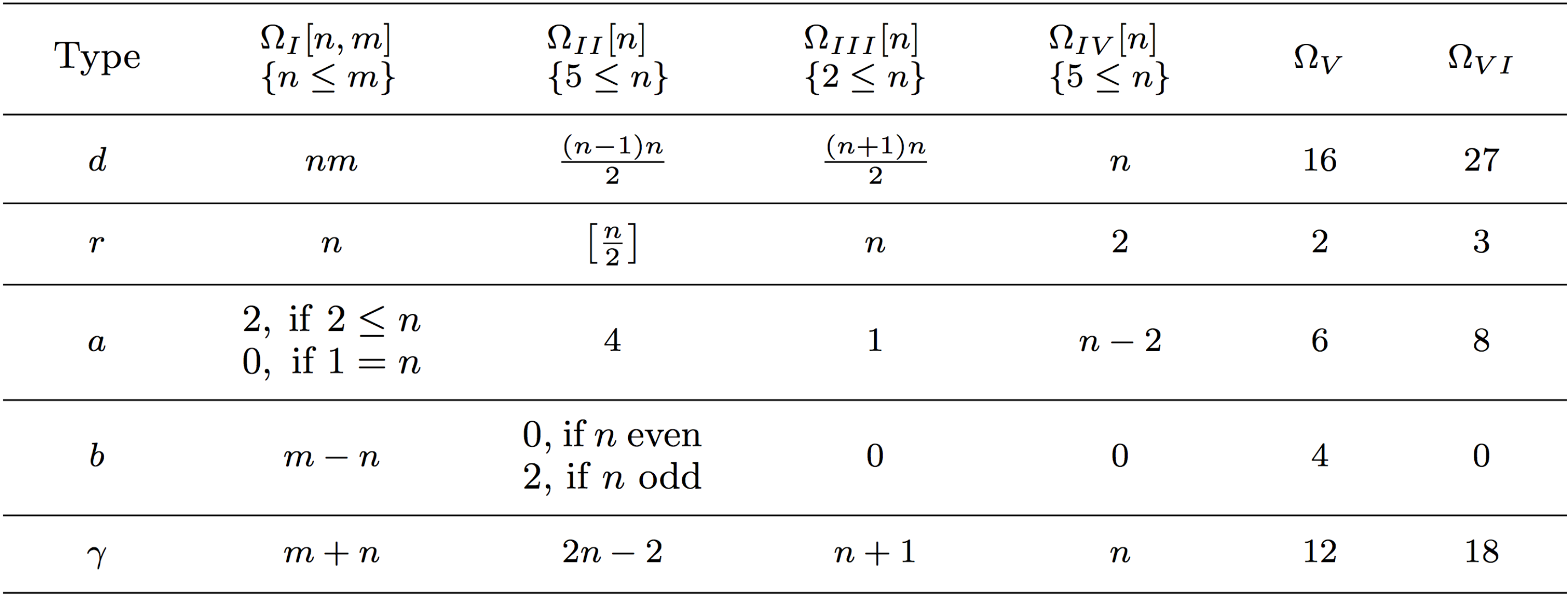}
\caption{Numerical invariants of irreducible  bounded symmetric domains.}
\end{center}
\end{figure}
\end{ex}

{\bf Acknowledgements.}
The author would like to thank Professor Andrea Loi for his interest in my work, his comments and various stimulating discussions.

\section{Proofs of Theorem \ref{prop ent bal}, \ref{thm dentropy000} and Corollary \ref{cor1}}

\subsection{Proof of Theorem \ref{prop ent bal}.}\label{proofthm1}
Assume that $\lmb\, g$ is balanced, namely that the $\ep$-function is constant. Consider  the \K\ potential $\lmb\, \phi \left( z\right) $ of $\lmb\, g$, let $\lmb\, \phi \left( z,\,\ov w\right) $ be its analytic continuation and let $K_{\lmb\,\phi}\left( z,\, \ov w\right) $ be  the reproducing kernel for the weighted Hilbert space $\hilb_{\lmb\, \phi}$. We have
\begin{equation}\label{epsilon1}
\ep_{\lmb\, g}\left( z,\, \ov w\right) =e^{-\lmb\,\phi\left( z,\, \ov w\right) }K_{\lmb\,\phi}\left( z,\, \ov w\right) =C>0,
\end{equation}
so $K_{\lmb\,\phi}\left( z,\, \ov w\right)$ never vanish. Hence, fixed a point $z_0$, the function
\begin{equation}\label{PSI00}
\psi\left( z,\, \ov w\right)=\frac{K_{\lmb\,\phi}\left( z,\,\ov w\right) K_{\lmb\,\phi}\left( z_0,\,\ov z_0\right)  }{K_{\lmb\,\phi}\left( z,\,\ov z_0\right) K_{\lmb\,\phi}\left( z_0,\,\ov w\right) }.
\end{equation}
is well defined. Note that $\psi\left( z_0,\,\ov w\right) =\psi\left( z,\,\ov z_0\right) =1$ and that
\begin{equation}\label{DIAST00}
\D_{z_0}\left( z\right) =\frac{1}{\lmb}\log\psi\left( z,\,\ov z\right) 
\end{equation}
is  the diastasis centered in $z_0$ associated to $g$. So the diastasis $\D_{z_0}$ is globally defined and $\lmb\,\D_{z_0}$ is a global \K\ potential for $\lmb\, g$. Consider the Hilbert space $\hilb_{\lmb\,\D_{z_0}}$ 
 given by:
\begin{equation*}
\hilb_{\lmb\,\D_{z_0}}=\left\{ f\in\ol\left( \Omega\right)  : \int_\Omega e^{-\lmb\,\D_{z_0}}\left|f\right|^2\, \frac{\omega^n}{n!}<\infty\right\}.
\end{equation*}
Let $K_{\lmb\,\D_{z_0}}\left( z,\, \ov w\right) $ be its reproducing kernel, since  the $\ep$-function does not depend on the \K\ potential, by \eqref{epsilon1} we have
\begin{equation}\label{epsilon2}
\ep_{\lmb\, g}\left( z,\, \ov w\right)=e^{-\lmb\,\D_{z_0}\left( z,\, \ov w\right)}K_{\lmb\,\D_{z_0}}\left( z,\, \ov w\right)=C
\end{equation}
where $\D_{z_0}\left( z,\, \ov w\right)$ denote the analytic continuation of $\D_{z_0}\left( z\right) =\D_{z_0}\left( z,\, \ov z\right) $.
In particular
\begin{equation*}
C= e^{-\lmb\,\D_{z_0}\left( z,\, \ov z_0\right)} K_{\lmb\,\D_{z_0}}\left( z,\, \ov z_0\right)=K_{\lmb\,\D_{z_0}}\left( z,\, \ov  z_0\right)  \in \hilb_{\lmb\,\D_{z_0}},
\end{equation*}
indeed by \eqref{PSI00} and \eqref{DIAST00} we have $\D_{z_0}\left( z,\, \ov z_0\right)= \frac{1}{\lmb} \log\psi\left( z,\,\ov z_0\right)=0$. It would be worth to say that $K_{\lmb\,\D_{z_0}} \in \hilb_{\lmb\,\D_{z_0}}$ by its very definition (indeed $K_{\lmb\,\D_{z_0}}\left( z,\, \ov  z_0\right)=\sum_{j=0}^{N} s_j\left( z\right)\, \ov {s_j\left( z_0\right)}$ where $\lbrace s_j \rbrace$ is a basis for $\hilb_{\lmb\,\D_{z_0}}$).
Thus $\hilb_{\lmb\,\D_{z_0}}$ contains the constant functions and
\[
\int_\Omega e^{-\lmb\,\D_{z_0}}\ \frac{\omega^n}{n!}<\infty.
\]
We conclude, by the very definition of the diastatic entropy, that $\Ent \left( \Omega,\, g\right) \left( z_0\right) <\lmb$.  Moreover, by \eqref{DIAST00}, the diastasis $\D_{z_0}$ is globally defined for any fixed $z_0 \in \Omega$ and  \cite[Proposition 2.2]{mossa1} says us that $\Ent \left( \Omega,\, g\right)$ is constant. As wished. 
\endproof

\subsection{Proof of Theorem \ref{thm dentropy000}.}
Recall that a homogeneous bounded  domain $\left( \Omega,\, g\right) $ is a bounded domain of $\C^n$ with a \K\ metric $g$ such that the group $G=\Aut (\Omega)\, \cap\, \Isom \left( \Omega,\, g\right) $ acts transitively on it, where $\Aut \left( \Omega\right) $ denotes the group of biholomorphisms of $\Omega$ and $\Isom \left( \Omega,\, g\right) $  the group of isometries of 
$\left( \Omega,\, g\right) $.
It is well-known that such a domain $\Omega$ is contractible and that $\omega=\frac{i}{2}\de \deb \phi$
for a globally defined \K\ potential $\phi$.  Indeed $\Omega$ is pseudoconvex being biholomorphically equivalent to a Siegel domain (see, e.g. \cite{vinberg} for a proof)  and so the existence of a global potential follows by a classical result of Hormander (see \cite{bulletin}) asserting that the equation $\deb u = f$ with $f$ $\deb$-closed form has a global solution on pseudoconvex domains
(see also \cite{mossa3, mossa2} and the proof of Theorem 4 in \cite{dlh}, for an explicit construction of the potential $\phi$ following  the ideas developed in \cite{dorf}).

We have to prove that if $1 > \Ent\left( \Omega,\, g\right) \left( z \right)$ then $g$ is balanced (the converse follows by Theorem 1). Assume that $1 > \Ent\left( \Omega,\, g\right) \left( z \right)$
i.e.
\[
\int_\Omega e^{-\D_{z}}\ \frac{\omega^n}{n!}<\infty,
\]
then the Hilbert space $\hilb_{\D_{z}}$ contains the constant functions and the associated $\ep$-function is strictly positive. By the homogeneity of $\left( \Omega,\, g\right) $ we conclude that $\ep_g$ has to be a positive constant (see the proof  of \cite[Lemma 1]{lm1}), that is $g$ is balanced. 
This concludes the proof of Theorem \ref{thm dentropy000}.
\subsection{Proof of Corollary \ref{cor1} }
In order to prove Corollary \ref{cor1} we need of the following result: 
\begin{lem}\label{ultimo lemma}(\cite[Theorem 2]{lm1} and \cite[(12)]{lm1}).
Let $\left( \Omega, \, \lmb\,g\right)$ be a homogeneous bounded domain and let $p_k,$ $b_k,$ $q_k$ and $\gamma_k$ be the constants defined by \eqref{def const} and \eqref{def const2} in Appendix. Then the metric $\lmb\, g$ is balanced if and only if 
\begin{equation*} 
 \lmb  >  \frac{1 + p_k + b_k + q_k/2}{\gamma_k}, \qquad 1 \le k \le r.
\end{equation*}
\end{lem}
By the very definition of diastatic entropy, we have that
\[
\Ent\left( \Omega, \, \lmb\,g\right) = \frac{\Ent\left( \Omega, \,g\right)}{\lmb},
\]
combining this with Theorem \ref{thm dentropy000} we obtain that $\Ent \left( \Omega, \, g \right) < \lmb$
 if and only if the metric $\lmb\, g$ is balanced.  Therefore 
 \[
\Ent\left( \Omega, \,g\right)=  \min \left\{ \lmb >0 \ | \ \lmb\, g  \text{ is balanced}\right\}, \qquad \forall\, z \in \Omega. \]
Together with Lemma \ref{ultimo lemma}, this concludes the proof of Corollary \ref{cor1}.

\section{Appendix}
In this section we give the definition of the constants
$p_k,$ $b_k,$ $q_k$ and $\gamma_k$ associated to the Piatetski-Shapiro's root structure (see \cite{PS69} for details).

Let $\left( \Omega \subset \C^n,\, g\right) $ be an homogenous bounded domain, and set $G=\Aut (\Omega) \cap \Isom  \left( \Omega ,\, g\right)$,
where $\Aut \left( \Omega\right) $ (resp. $\Isom \left( \Omega ,\, g\right)$) denotes the group of invertible holomorphic maps (resp. $g$-isometries) of $\Omega$.
By \cite[Theorem 2 (c)]{dorf},
 there exists a connected split solvable Lie subgroup
 $S \subset G$
 acting simply transitively on the domain $\Omega$. Fix a point $p_0 \in \Omega$,
 we have a diffeomorphism $\iota: S \f \Omega$ with
$\iota \left( \g\right)  = \g \cdot p_0$. By  differentiation,
 we get a linear isomorphism between the Lie algebra $\gs$ of $S$ and the tangent space $T_{p_0}\HBD$, given by
 $X
  \overset{\sim}{\mapsto} X \cdot p_0 $.
  Then the evaluation of the K\"ahler form
 $\omega$ on $\gs$
 is given by
\begin{equation*} \label{eqn:beta}
 \omega\left( X\cdot p_0, Y \cdot p_0\right)  = \beta\left( \left[X,\,Y\right]\right) \qquad
 X, Y \in \gs
\end{equation*}
 with a certain linear form $\beta \in \gs^*$.
Let $j : \gs \to \gs$
 be the linear map defined in such a way that
 \[\left( j\,X\right)  \cdot p_0 = \sqrt{-1} \left( X \cdot p_0\right)  \text{ for } X \in \gs.\]
We have
 $g\left( X \cdot p_0,\,Y \cdot p_0\right)  = \beta\left( \left[jX,\, Y\right]\right) $ 
 for $X,\, Y \in \gs$,
 so the right-hand side defines a positive inner product on $\gs$.
Let $\mathfrak{a}$ be the orthogonal complement of $\left[\gs,\, \gs\right]$ in $\gs$
 with respect to this inner product.
Then $\mathfrak{a}$ is a commutative Cartan subalgebra of $\gs$.
Given a linear form $\alpha$ on the Cartan algebra $\mathfrak{a}$,
 we denote by $\gs_{\alpha}$ the subspace
 \[\gs_{\alpha} = \set{X \in \gs}{\left[C,\,X\right] = \alpha\left( C\right) X, \ \forall\, C \in \mathfrak{a}}\]
 of $\gs$.
We say that $\alpha$ is a \textit{root}
 if $\gs_{\alpha} \ne \{0\}$ and $\alpha \ne 0$.
Thanks to \cite[Chapter 2, Section 3]{PS69} or \cite[Theorem 4.3]{RV73},
 there exists a basis
 $\{\alpha_1, \dots, \alpha_r\}$, $(r := \dim \mathfrak{a})$
 of $\mathfrak{a}^*$
 such that every root is one of the following:
\[
 \alpha_k, \ \alpha_k/2, \ (k=1, \dots, r), \qquad (\alpha_l \pm \alpha_k)/2, \qquad 1 \le k < l \le r.
\]
For $k=1, \dots, r$ we define (see \cite[Definition 4.7]{RV73} and \cite[(2.7)]{N01})
\begin{equation}\label{def const}
 p_k := \sum_{i<k} \dim \gs_{\left( \alpha_k - \alpha_i\right) /2}, \quad
 q_k := \sum_{l>k} \dim \gs_{\left( \alpha_l - \alpha_k\right) /2}, \quad
 b_k := \frac{1}{2}\dim \gs_{\alpha_k/2}, 
\end{equation}
\begin{equation}\label{def const2}
\gamma_k := 4\, \beta\left( \left[jA_k,\, A_k\right]\right),
\end{equation}
where $\lbrace A_1, \dots , A_r \rbrace$ is the basis of $\mathfrak{a}$ dual to $\lbrace \alpha_1, \dots , \alpha_r \rbrace$.

\end{document}